\documentclass[11pt]{article}
\usepackage[latin1]{inputenc}
\usepackage{amssymb}
\usepackage{amsmath}
\usepackage{amsthm}

\input pictex.tex
\newcommand{\R}{\mathbb{R}}

\newcommand{\N}{\mathbb N}
\newcommand{\Z}{\mathbb Z}

\theoremstyle{definition}

\newtheorem{thm}{Theorem}[section]

\newtheorem{prop}[thm]{Proposition}
\newtheorem{defn}[thm]{Definition}

\newtheorem{lem}[thm]{Lemma}

\newtheorem{rem}[thm]{Remark}

\newtheorem{ex}[thm]{Example}

\newcommand{\vs}{\vspace{0.3cm}}

\newcommand{\vsp}{\vspace{0.1cm}}

\date{}
\author{}

\begin{document}

\title{Left-orderings on free products of groups}
\author{Crist\'obal Rivas}
\maketitle
\begin{abstract} \noindent We show that no left-ordering on a free
product of (left-orderable) groups is isolated. In particular, we
show that the space of left-orderings of a free product of finitely
generated groups is homeomorphic to the Cantor set. With the same
techniques, we also give a new and constructive proof of the fact
that the natural conjugation action of the free group (on two or
more generators) on its space of left-orderings has a dense orbit.
\end{abstract}

%%%%%%%%%%%%%%%%%%%%%%%%%%%%%%%%%%%%%%%%%%%%%%%%%%%%%%%%%%%%%%%%%%%%%%%
\vs

\noindent{\Large \bf Introduction}

\vs
%%%%%%%%%%%%%%%%%%%%%%%%%%%%%%%%%%%%%%%%%%%%%%%%%%%%%%%%%%%%%%%%%%%%%%%%

A (non-necessarily total) order relation $\preceq$ on a group
$\Gamma$ is said to be a {\em partial-left-ordering}, if for every
$\gamma_1, \gamma_2,\gamma_3$ in $\Gamma$, we have that
$\gamma_1\prec \gamma_2$ implies $\gamma_3\gamma_1\prec
\gamma_3\gamma_2$. An element $\gamma \in \Gamma$ is called
$\preceq$-positive (resp. $\preceq$-negative) if $id\prec \gamma$
(resp. $\gamma \prec id$). The subset of $\preceq$-positive
elements, usually called the {\em positive cone} for $\preceq$, will
be denoted by $P_\preceq$. Clearly, $P_\preceq$ satisfies

\vsp\vsp

\noindent $(O1)$ $P_\preceq  P_\preceq \subseteq P_\preceq \;$, that
is, $P_\preceq$ is a semi-group, and

\vsp

\noindent $(O2)$ $P_\preceq \cap P_\preceq^{-1}=\emptyset$, where
$P^{-1}_\preceq=\{g^{-1}\in \Gamma\mid g\in P_\preceq\}=\{g\in
\Gamma\mid g\prec id\}$.

\vsp\vsp

If in addition, $\preceq$ is a total order, we will simply say that
$\preceq$ is a {\em left-ordering}. In this case, the set of
$\preceq$-positive elements also satisfies

\vsp\vsp

\noindent $(O3)$ $\Gamma=P_\preceq \cup P^{-1}_\preceq \cup \{id\}$.

\vsp\vsp

Conversely, given any subset $P\subseteq \Gamma$ satisfying the
conditions $(O1)$, $(O2)$ and $(O3)$ (resp. $(O1)$ and $(O2)$)
above, we can define a left-ordering (resp. a partial-left-ordering)
$\preceq_P$ by letting $f\prec_P g$ if and only if $f^{-1}g\in P$.
We will usually identify $\preceq$ with $P_{\preceq}$.

\vsp

Given a group $\Gamma$ (of arbitrary cardinality), we denote the set
of all partial-left-orderings on $\Gamma$ by
$\mathcal{PLO}(\Gamma)$. This set has a natural topology first
exploited by Sikora for the case of (total orderings on) countable
groups \cite{sikora}. This topology can be defined by identifying
$P\in \mathcal{PLO}(\Gamma)$ with its characteristic function
$\chi_P \in \{0,1\}^\Gamma $. In this way, we can view
$\mathcal{PLO}(\Gamma)$ embedded in $\{0,1\}^\Gamma$. This latter
space, with the product topology, is a Hausdorff, totally
disconnected, and compact space. It is not hard to see that (the
image of) $\mathcal{PLO}(\Gamma)$ is closed inside, and hence
compact as well (see \cite{linnell,witte,navas,sikora} for details).
In the same way, for a left-orderable group $\Gamma$, the space of
all left-orderings, here denoted $\mathcal{LO}(\Gamma)$, is closed
inside $\mathcal{PLO}(\Gamma)$, hence compact as well. In
\cite{linnell}, it is shown that $\mathcal{LO}(\Gamma)$ is either
finite or uncountable.

\vsp

A basis of neighborhoods of $\,\preceq \,$ in $\,
\mathcal{LO}(\Gamma)$ is the family of the sets $\,
V_{f_1,\ldots,f_k}=\{ \preceq^\prime \in \mathcal{LO}(\Gamma) \,\mid
\, id \prec^\prime f_i, \text{ for } i=1,\ldots k\} \,$, where
$\{f_1,\ldots,f_k\}$ runs over all finite subsets of
$\,\preceq$-positive elements of $\Gamma$ (the same being true for
$\mathcal{PLO}(\Gamma)$). Therefore, it is natural to say that a
left-ordering $\preceq$ of $\Gamma$ is {\em isolated} if and only if
there is a finite family $\{\gamma_1,\ldots,\gamma_n\}\subset
\Gamma$ such that $\preceq$ is the only left-ordering of $\Gamma$
with the property that $\gamma_i\succ id$, for $1\leq i \leq n$.

\vsp

Knowing whether a given group has an isolated left-ordering turns
out to be a natural and old question in the theory of left-orderable
groups (although not always expressed in topological terms...). A
major progress in the understanding of groups having isolated
left-orderings, is the classification of groups admitting only
finitely many left-orderings (all of them isolated) made by Tararin
\cite[Theorem 5.2.1]{kopytov}. In addition, we count with the
remarkable examples of groups admitting infinitely many
left-orderings together with some isolated left-orderings, such as
the braid groups \cite{dd} (see however \cite{navas bertol}), and
the groups appearing in \cite{ito,navas-hecke}. On the other hand,
it is known that some classes of groups, such as nilpotent groups
\cite{navas} (more generally, left-orderable groups of
sub-exponential growth \cite[Remark 2.2.3]{tesis}) and the groups
appearing in \cite{rivas2}, have no isolated left-orderings unless
they have only finitely many left-orderings.

\vsp

In the case of the free group of finite rank $F_n$, $n\geq2$, it was
proved by McCleary \cite{mccleary} that $F_n$ has no isolated
left-orderings\footnote{The fact that the free groups of infinite
rank has no isolated left-orderings is easy and appears, for
instance, in \cite{braids}.}. McCleary's proof relies on the study
of the so called free-lattice-ordered group (in his case) over the
free group, which is a universal object introduced by Conrad in
\cite{conrad2}. An independent proof of this fact was given by Navas
in \cite{navas}, where he studies the so-called {\em dynamical
realization of a left-ordering} (see \S \ref{dinamical realization})
of $F_n$, which is an order-preserving action on the real line that
encodes all the information of the given left-ordering.

\vsp

In this article, we  simplify and generalize Navas' approach to get
a generalization of McCleary's result for the case of free products
of left-orderable groups. (Recall that the free product of
left-orderable groups is left-orderable \cite[Corollary
6.1.3]{kopytov}.) We show

\vsp\vsp

\noindent {\bf Theorem A:} {\em Let $G$ and $H$ be two
left-orderable groups. Then the free product $G*H$ has no isolated
left-orderings.}

\vsp\vsp

To prove Theorem A we first work the case where $G$ and $H$ are
finitely generated \S \ref{sec finito generado}. Then, in \S
\ref{sec extension por compacidad}, we use the compactness of
$\mathcal{PLO}(\Gamma)$ to provide an argument ensuring Theorem A.
We note that Theorem A does not extends to the case of amalgamated
free products, since the groups with isolated left-orderings
appearing in \cite{ito,navas-hecke} (for instance, the braid group
$B_3$) are of that form.

\vsp

A direct consequence of Theorem A is that no positive cone of a
left-ordering on a free product of groups is finitely generated as a
semigroup (see for instance \cite[Proposition 1.8]{navas}). However,
the converse to this is not true. In \S \ref{sec ejemplo}, we show
that $\langle a, b \mid bab^{-1}=a^{-2}\rangle$ is a group with an
isolated left-ordering  whose positive cone is not finitely
generated as a semigroup.

\vsp

Besides its compactness, $\mathcal{LO}(\Gamma)$ has another very
important property, namely, that the group $\Gamma$ naturally acts
on it by conjugation:
$$\gamma(\preceq)=\preceq_\gamma,\text{ where }\gamma_1\prec_\gamma
\gamma_2 \text{ if and only if } \gamma \gamma_1 \gamma^{-1}\prec
\gamma \gamma_2 \gamma^{-1}.$$ This action turns out to be by
homeomorphisms since
$\gamma(V_{\gamma_1,\ldots,\gamma_k})=V_{\gamma_1\gamma^{-1},\ldots,
\gamma_k\gamma^{-1}}$. This action was defined by Ghys and was first
exploited in \cite{witte} by Morris-Witte.

\vsp

In \cite{clay 2}, Clay found a strong connection between the
conjugation action of $\Gamma$ on its space of left-orderings and
some natural representations of the free-lattice-ordered group over
$\Gamma$. In the special case of a free group, this connection,
together with a previous result of Kopytov \cite{kopytov2}, allowed
him to show

\vsp\vsp

\noindent {\bf Theorem B (Clay):} {\em Let $\mathcal{F}$ be  a free
group of countable rank greater than one. Then, the space of
left-orderings of $\mathcal{F}$ has a dense orbit under the natural
conjugation action of $\mathcal{F}$.}

\vsp\vsp

Nevertheless, his proof is highly non-constructive, and Kopytov's
result also involves the free-lattice-ordered group over the free
group. In Section \ref{sec dense} of this work, we use our dynamical
machinery to give an explicit and self-contained construction of a
left-ordering on $\mathcal{F}$ whose set of conjugates is dense.
However, our method does not solve the following question, that may
have some interest in rigidity theory.

\vsp

\noindent{\bf Question:} Does $\mathcal{F}$ admits a dense orbit for
the diagonal action on $\mathcal{LO}(\mathcal{F})\times
\mathcal{LO}(\mathcal{F})$?

\vs

\noindent {\bf Acknowledgments:}

\vs

I would like to thanks Luis Paris for his invitation to the
Université de Bourgogne where the finitely generated version of
Theorem A was found, and Christian Bonnatti for his interest on
left-orderable groups and for the valuable comments that makes this
article possible. I am also grateful to Andrés Navas for explaining
me how to pass from the finitely generated version of Theorem A to
the final version of it, and to Adam Clay for his interest on this
subject and some corrections on an earlier draft of this work.

%%%%%%%%%%%%%%%%%%%%%%%%%%%%%%%%%%%%%%%%%%%%%

\section{The dynamical realization of a
left-ordering} \label{dinamical realization}

%%%%%%%%%%%%%%%%%%%%%%%%%%%%%%%%%%%%%%%%%%%%%%%%%%%%%%%%%%%%%%

\hspace{0.35 cm} Though orderability may look as a very algebraic
concept, it has a deep (one-dimensional) dynamical content. For
instance, a group is left-orderable if and only if it embeds in the
group of order-preserving automorphisms of a totally ordered set
$\Omega$; see for instance \cite[Theorem 3.4.1]{kopytov}.

\vsp

For the case of countable groups ({\em e.g.} finitely generated), we
can give more dynamical information since we can take  $\Omega$ as
being the real line (see \cite[Theorem 6.8]{ghys}, or \cite{navas}
for further details).

\vsp

\begin{prop} \label{real din}{\em For a countable infinite group $\Gamma$, the following two properties are
equivalent:\\

\noindent -- $\Gamma$ is left-orderable,\\

\noindent -- $\Gamma$ acts faithfully on the real line by
orientation-preserving homeomorphisms. That is, there is an
homomorphic embedding $\Gamma\to Homeo_+(\R)$.}
\end{prop}

\noindent\textit{Sketch of proof: } To show that a subgroup of
$Homeo_+(\R)$ is left-orderable, we construct what is usually called
an {\em induced left-ordering}. To do this, we take a dense sequence
$(x_0,x_1,\ldots)$ of points in $\R$, and we define
$\preceq_{(x_0,x_1,\ldots)}$ by declaring
$$\gamma\succ_{(x_0,x_1,\ldots)} id\;\; \text{ if and only if }
\;\;\gamma(x_i)>x_i,$$ where $i=min\{ j\mid x_j\not=\gamma(x_j)\}$.
Showing that $\preceq_{(x_0,x_1,\ldots)}$ is a total left-ordering
is routine.

\vsp

For the converse, we construct what is called \textit{a dynamical
realization of a left-ordering $\preceq$.} Fix an enumeration
$(\gamma_i)_{i \geq 0}$ of $\Gamma$ such that $\gamma_0=id$, and let
$t_\preceq(\gamma_0)=0$. We shall define an order-preserving map
$t_\preceq: \Gamma \to \R$ by induction. Suppose that
$t_\preceq(\gamma_0),
t_\preceq(\gamma_1),\ldots,t_\preceq(\gamma_i)$ have been already
defined. Then if $\gamma_{i+1}$ is greater (resp. smaller) than all
$\gamma_0,\ldots, \gamma_i$, we define $t_\preceq(\gamma_{i+1})=
max\{t_\preceq(\gamma_0),\ldots, t_\preceq(\gamma_i)\}+1$ (resp.
$min\{t_\preceq(\gamma_0),\ldots, t_\preceq(\gamma_i)\}-1$). If
$\gamma_{i+1}$ is neither greater nor smaller than all
$\gamma_0,\ldots,\gamma_i$, then there are
$\gamma_n,\gamma_m\in\{\gamma_0,\ldots , \gamma_i \}$ such that
$\gamma_n\prec \gamma_{i+1}\prec \gamma_m$ and no $\gamma_j$ is
between $\gamma_n,\gamma_m$ for $0\leq j\leq i$. Then we set
$t_\preceq(\gamma_{i+1})=(t_\preceq(\gamma_n)+t_\preceq(\gamma_m))/2$.

\vsp

Note that $\Gamma$ acts naturally on $t_\preceq(\Gamma)$ by
$\gamma(t(\gamma_i)) = t_\preceq(\gamma\gamma_i)$, and that this
action extends continuously to the closure of $t_\preceq(\Gamma)$.
Finally, one can extend the action to the whole real line by
declaring the map $\gamma$ to be affine on each interval of the
complement of $\overline{t_\preceq(\Gamma)}$. $\hfill\square$

\vs

We have just constructed an embedding of a countable, left-ordeable
group $ \Gamma $ into $Homeo_+(\R)$. We call this embedding a
dynamical realization of the left-ordered group $(\Gamma,\preceq)$.
The order preserving map $t_\preceq$ is called the reference map.

\vsp

\begin{rem} \label{rem conj1} As constructed above, the dynamical realization depends
not only on the left-ordering $\preceq$, but also on the enumeration
$(\gamma_i)_{i\geq 0}$. Nevertheless, it is not hard to check that
dynamical realizations associated to different enumerations (but the
same ordering) are \textit{topologically conjugate}.\footnote{Two
actions $\phi_1\!: \Gamma \to \mathrm{Homeo}_+(\R)$ and
$\phi_2\!:\Gamma \to \mathrm{Homeo}_+(\R)$ are topologically
conjugate if there exists $\varphi \in \mathrm{Homeo}_+(\R)$ such
that $\varphi\circ \phi_1(\gamma) = \phi_2(\gamma) \circ \varphi$
for all $\gamma \in \Gamma$.} Thus, up to topological conjugacy, the
dynamical realization depends only on the ordering $\preceq$ of
$\Gamma$.

\vsp

An important property of dynamical realizations is that they do not
admit global fixed points (\textit{i.e.,} no point is stabilized by
the whole group). Another important property is that
$0=t_\preceq(id)$ has a {\em free orbit} (\textit{i.e} $\{\gamma\in
\Gamma\mid \gamma(t_\preceq(id))=t_\preceq(id)\}=\{id\}$ ). Hence
$\gamma\succ id$ if and only if $\gamma(t_\preceq(id))=\gamma(0)>0=
t_\preceq(id)$, which allows us to recover the left-ordering from
its dynamical realization.
\end{rem}

The following well-known Proposition will serve us to approximate a
given left-ordering by looking at its dynamical realization. For the
reader convenience, we sketch the proof below.

\begin{prop}\label{lema extension} {\em Let $\Gamma$ be a left-orderable group, and let
$D:\Gamma \to Homeo_+(\R)$ be a (not necessarily faithful)
homomorphism. Let $x_0\in \R$ and let $\preceq_{x_0}$ be the
partial-left-ordering defined by $\gamma \succ_{x_0} id$ if and only
if $D(\gamma)(x_0)>x_0$. Then $\preceq_{x_0}$ can be extended to a
(total) left-ordering $\preceq$ such that $\gamma\succ_{x_0} id$
implies $\gamma \succ id$.}
\end{prop}

\noindent {\em Sketch of proof:} Let $H=\{\gamma\in \Gamma \mid
D(\gamma)(x_0)=x_0\}$. Let $\preceq^\prime$ be any left-ordering on
$H$. Define $\preceq$ by
$$ g\succ id \Leftrightarrow \left\{
\begin{array}{l } D(\gamma)(x_0)>x_0  \text{ or } \\
D(\gamma)(x_0)=x_0 \;\;\text{ and } g\succ^\prime id .
\end{array} \right.$$
Showing that $\preceq $ is a left-ordering on $\Gamma$ is
straightforward. $\hfill\square$

\begin{defn} \label{def din rel like} Let $\preceq$ be a left-ordering on a countable group
$\Gamma$. Let $D:\Gamma \to Homeo_+(\R)$ be an homomorphic embedding
with the property that {\em there exists $x \in \R$ such that, for
$\gamma_1$ and $\gamma_2$ in $\Gamma$, we have that $\gamma_1\prec
\gamma_2$ if and only if $D(\gamma_1)(x)<D(\gamma_2)(x)$.} We call
$D$ a {\em dynamical realization-like homomorphism} for $\preceq$.
The point $x$ is called {\em reference point} for $D$.
\end{defn}

\begin{ex} The embedding given by any dynamical realization of any
countable left-ordered group $(\Gamma,\preceq)$ is a dynamical
realization-like homomorphism for $\preceq$ with reference point
$0=t_\preceq(id)$.
\end{ex}

\begin{rem}\label{rem cajas}
Note that, if $D$ is a dynamical realization-like homomorphism for
$\preceq$, with reference point $x$, and if $\varphi:\R\to \R$ is
any increasing homeomorphism, then the conjugated homomorphism
$D_\varphi$ defined by $D_\varphi(g)= \varphi D(g) \varphi^{-1}$ is
again a dynamical realization-like homomorphism for $\preceq$ but
with reference point $\varphi(x)$.
\end{rem}

For the rest of this section, $\Gamma$ will be a countable (not
necessarily finitely generated) left-orderable group, and $\Gamma_0$
a finite subset of $\Gamma$ such that $\Gamma_0=\Gamma_0^{-1}$. We
will also denote $\langle \Gamma_0\rangle$ the subgroup generated by
$\Gamma_0$. Finally, for $w\in \langle \Gamma_0 \rangle$, we will
denote by $|w|_{\Gamma_0}$ the word length of $w$ with respect to
$\Gamma_0$.

\vs

The following notion will be essential in our work.

\begin{defn} \label{def box} Let $B_{\Gamma_0}(n)=\{w\in \langle\Gamma_0\rangle \mid
|w|_{\Gamma_0}\leq n\}$ be the {\em ball} of radius $n$ in $\langle
\Gamma_0 \rangle$. Given $B_{\Gamma_0}(n)\subseteq \Gamma$ and a
left-ordering $\preceq$ of $\Gamma$, we let
$$\lambda_{(B_{\Gamma_0}(n),\preceq)}^-= \min_\preceq \{ w\in B_{\Gamma_0}(n)\} ,\;\;\;\; \lambda_{(B_{\Gamma_0}(n),\preceq)}^+= \max_\preceq \{ w\in B_{\Gamma_0}(n) \}.$$

Now, let $D$ be a dynamical realization-like homomorphism for
$\preceq$, with reference point $x$. Then, we will refer to the
square
$[D(\lambda_{(B_{\Gamma_0}(n),\preceq)}^-)(x),D(\lambda_{(B_{\Gamma_0}(n),\preceq)}^+)(x)]^2\subset
\R^2$ as the $(B_{\Gamma_0}(n),\preceq)$-box.
\end{defn}

\begin{rem} \label{rem delta} Note that, from the left-invariance of
$\preceq$, we have that
$\big|\lambda^\pm_{(B_{\Gamma_0}(n),\preceq)}\big|_{\Gamma_0}=n$,
and that there is $\delta^+_n\in \Gamma_0$ (resp. $\delta_n^-\in
\Gamma_0$) such that $\delta^+_n
\lambda^+_{(B_{\Gamma_0}(n),\preceq)}=\lambda^+_{(B_{\Gamma_0}(n+1),\preceq)}$
(resp. $\delta^-_n
\lambda^-_{(B_{\Gamma_0}(n),\preceq)}=\lambda^-_{(B_{\Gamma_0}(n+1),\preceq)}$).

\end{rem}

\vsp

Now let $\preceq$ be a left-ordering on $\Gamma$. The next lemma
shows that the $(B_{\Gamma_0}(n),\preceq)$-box contains the
information of the $\preceq$-signs of the elements in
$B_{\Gamma_0}(n)$.

\begin{lem} \label{signo}   \label{remark signo}
{\em Let $D:\Gamma\to Homeo_+(\R)$ be a dynamical realization-like
homomorphism for $\preceq$ with reference point $x$. Then, for every
$w_1$ and $w_2$ in $B_{\Gamma_0}(n)$, we have that $D(w_1)(x)$
belongs to
$[D(\lambda^-_{(B_{\Gamma_0}(n),\preceq)})(x),D(\lambda^+_{(B_{\Gamma_0}(n),\preceq)})(x)]$,
and $D(w_1)(x)>D(w_2)(x)$ if and only $w_1 \succ w_2$.

\vsp

Moreover, for any representation $\tilde{D}:\Gamma\to Homeo_+(\R)$
such that, for every $\gamma\in \Gamma_0$, the  graphs\footnote{As
usual, for $f\in Homeo_+(\R)$, the set $\{(x,f(x))\mid x\in \R\}
\subset \R^2$ is called the graph of $f$.} of $\tilde{D}(\gamma)$
coincide with the graphs of $D(\gamma)$ inside
$[D(\lambda_{(B_{\Gamma_0}(n),\preceq)}^-)(x),D(\lambda_{(B_{\Gamma_0}(n),\preceq)}^+)(x)]^2$,
we have that $D(w)(x)=\tilde{D}(w)(x)$ for all $w\in
B_{\Gamma_0}(n)$.}
\end{lem}

\noindent {\em Proof:} From Definition \ref{def din rel like}, it
follows that for any $w_1$ and $w_2$ in $\Gamma$,
$D(w_1)(x)>D(w_2)(x)$ if and only if $w_1\succ w_2$. Now, for $w\in
B_{\Gamma_0}(n)$, we have that
$\lambda^{-}_{(B_{\Gamma_0}(n),\preceq)} \preceq w\preceq \lambda^{+
}_{(B_{\Gamma_0}(n),\preceq)}$. In particular, $D(w)(x)\in
[D(\lambda^-_{(B_n,\preceq)})(x),D(\lambda^+_{(B_n,\preceq)})(x)] $,
which shows the first part of the lemma.

\vsp

To show the second part, we note that every initial segment $w_1$ of
any reduced\footnote{By ``reduced" we mean a word of minimal length
among words in $\Gamma_0$.} word $w\in B_{\Gamma_0}(n)$ lies again
in $B_{\Gamma_0}(n)$. Hence, if $w=\alpha_j\ldots \alpha_1$, $j\leq
n$, where $\alpha_i\in \Gamma_0=\Gamma_0^{-1}$, is a reduced word,
then the points $x_1=D(\alpha_1)(x)$, $\;x_2=D(\alpha_2)(x_1)$,
\ldots , $x_j=D(\alpha_j)(x_{j-1})=D(w)(x)$, they all belong to
$[D(\lambda^-_{(B_n,\preceq)})(x),D(\lambda^+_{(B_n,\preceq)})(x)]$.
In particular, $x_1=D(\alpha_1)(x)=\tilde{D}(\alpha_1)(x), \ldots ,
x_j=D(\alpha_j)(x_{j-1})=\tilde{D}(\alpha_j)(x_{j-1})$, which shows
that $D(w)(x)=\tilde{D}(w)(x)$. $\hfill\square$

%%%%%%%%%%%%%%%%%%%%%%%%%%%%%%%%%%%%%%%%%%%%%%%%%%%%%%%%%%%%%%%%%%%%%%%%%

\section{Proof of Theorem A}

\subsection{The case where $\Gamma=G*H$ is finitely generated}
\label{sec finito generado}

%%%%%%%%%%%%%%%%%%%%%%%%%%%%%%%%%%%%%%%%%%%%%%%%%%%%%%%%%%%%%%%%%%%%%%%%%%%

\hspace{0.4 cm} Recall that the space of left-orderings of a
countable group $\Gamma$ is metrizable \cite{witte,navas,sikora}.
For instance, if $\Gamma$ is finitely generated, and $B_n$ denote
the ball of radius $n$ with respect to a finite generating set, then
we can declare $dist(\preceq_1,\preceq_2)=1/n$, if $B_n$ is the
largest ball on which $\preceq_1$ and $\preceq_2$ coincide. In
particular, if $\mathcal{LO}(\Gamma)$ contains no isolated points,
then $(\mathcal{LO}(\Gamma),dist)$ becomes a compact, Hausdorff and
locally disconnected metric space that has no isolated points. Hence
it is homeomorphic to the Cantor set \cite{hockin young}.

\vsp

For the rest of this section, $\Gamma$ will be the free product
$G*H$. Both groups $G$ and $H$ are assumed to be finitely generated
and left-orderable. The generating set of $G$ and $H$ will be
denoted $G_0=\{g_1, \ldots, g_k\}$  and $H_0=\{h_1,\ldots h_\ell\}$
respectively. We assume that $G_0$ and $H_0$ are closed under
inversion. In particular, $\Gamma=G*H$ is generated by
$\Gamma_0=\{g_1, \ldots, g_k,h_1,\ldots h_\ell\}=\Gamma_0^{-1}$.
Since in this case we have that $\langle \Gamma_0\rangle=\Gamma$, we
will denote the sets $B_{\Gamma_0}(n)$ (see Definition \ref{def
box}) simply by $B_n$.

\begin{thm} \label{teo finito generado}{\em No left-ordering on $G*H$ is isolated. In particular,
$\mathcal{LO}(G*H)$ is homeomorphic to the Cantor set.}
\end{thm}

\noindent {\em Proof:} To prove Theorem \ref{teo finito generado},
it is enough to show that, given a left-ordering $\preceq$ and a
finite subset $F$ of  $\Gamma$, there is a left-ordering
$\preceq^\prime$ different from $\preceq$ such that $\preceq^\prime$
coincides with $\preceq$ over $F$.

\vsp

To show this, we will perform a perturbation of the dynamical
realization $D:\Gamma\to Homeo_+(\R)$ of $\preceq$. This
perturbation will be made by conjugating the action of one of the
factors by an order preserving homeomorphism $\varphi:\R \to \R$,
while keeping the action of the second factor untouch. As explained
in Remak \ref{rem conj1}, we have that $\gamma\succ \gamma^\prime$
if and only if $D(\gamma)(0)>D(\gamma^\prime)(0)$ for all $\gamma,
\gamma^\prime$ in $\Gamma$.

\vsp

Since $\preceq$ is fixed, to avoid heavy notation, we will denote
the elements $\lambda^\pm_{(B_n,\preceq)}$ simply by
$\lambda^\pm_n$. We now let $n\in \N$ be such that $F\subseteq B_n$.

\vsp

Now, consider $\lambda^+_{n+1}$, and let $g\in G_0$ and $h\in H_0$
be such that $g\lambda^+_{n+1}\succ\lambda^+_{n+1}$ and
$h\lambda^+_{n+1}\succ\lambda^+_{n+1}$. Since we are not making any
different assumption on $G$ and $H$, we can assume that
$g\lambda^+_{n+1}\succ h\lambda^+_{n+1}$ (otherwise we change the
names...).

\vsp

We also let $x_0,x_1, y_0, y_1$ in $\R$ be such that
$$D(\lambda_{n+1}^+)(0)<x_0<x_1<D(h\lambda_{n+1}^+)(0) <
D(g\lambda_{n+1}^+)(0)<y_1<y_0.$$

\vsp

We let $\varphi\in Homeo_+(\R)$ be such that $supp(\varphi)=\{x\in
\R\mid \varphi(x)\not=x\}=(x_0,y_0)$ and $\varphi(x_1)>y_1$. This
implies that
\begin{equation}\label{distintos}\varphi \circ D(h\lambda_{n+1}^+) \circ
\varphi^{-1}(0)> D(g\lambda_{n+1}^+)(0),\end{equation} where $\circ$
is the composition operation. Moreover, for any $\bar{h}\in H_0$ and
any $x\in [D(\lambda_n^-)(0),D(\lambda_n^+)(0)]$, we have that
$D(\bar{h})(x)\leq D(\lambda_{n+1}^+)(0)<x_0$. Thus we conclude,
\begin{equation}\varphi\circ D(\bar{h})\circ \varphi^{-1}(x)=D(\bar{h})(x),
 \text{ for all } x\leq D(\lambda_n^+)(0) \text{ and all } \bar{h}\in H_0.\label{critico}\end{equation}

\vsp

Now, let $D_\varphi:\Gamma \to Homeo_+(\R)$ be defined by
$D_\varphi(\bar{g})=D(\bar{g})$ for all $\bar{g}\in G$, and
$D_\varphi(\bar{h})=\varphi\circ D(\bar{h}) \circ \varphi^{-1}$ for
all $\bar{h}\in H$. Since $\Gamma$ is the free product of $G$ and
$H$, we have that $D_\varphi$ is an homomorphism (not necessarily
injective). Now, from the definition of $D_\varphi$ and equation
(\ref{critico}), we have that
\begin{equation}\label{la conju} D(\gamma)(x)=D_\varphi(\gamma)(x) \text{ for any $\gamma \in
\Gamma_0$ and any $x\leq D(\lambda_n^+)(0)$.}\end{equation}  In
particular, for each $\gamma\in \Gamma_0$, the graphs of $D(\gamma)$
and $D_\varphi(\gamma)$ coincide inside the square \break
$[D(\lambda_n^-)(0),D(\lambda_n^+)(0)]^2$. Hence, from Lemma
\ref{signo}, we conclude that
\begin{equation} \label{cerca} \text{ for all } \gamma\in B_n,
\;\;D(\gamma)(0)=D_\varphi(\gamma)(0).\end{equation}

\vsp

Now, from Lemma \ref{lema extension}, we have that there is a
left-ordering $\preceq^\prime$ on $\Gamma$ such that
$D_\varphi(\gamma)(0)>0$ implies $\gamma\succ^\prime id$. Then,
equation (\ref{cerca}) implies that $\preceq$ and $\preceq^\prime$
coincide on $B_n$, hence, on $F$.

\vsp

However, if we let $\delta_n \in \Gamma_0$ be such that
$\delta_n\lambda^+_n=\lambda_{n+1}^+$ (see Remark \ref{rem delta}),
we have that $D(g\lambda_{n+1})(0)=D(g)\circ D(\delta_n)\circ
D(\lambda_n^+)(0)$. Hence, from the definition of $D_\varphi$ and
equations (\ref{la conju}) and (\ref{cerca}), we conclude that
$D(g\lambda_{n+1})(0)=D_\varphi(g\lambda_{n+1})(0)$. Moreover, from
the definition of $\varphi$, we have that $\varphi \circ
D(h\lambda_{n+1}^+) \circ \varphi^{-1}(0) = D_\varphi(h)\circ
\varphi  \circ D(\lambda_{n+1}^+) \circ\varphi^{-1}(0)=D_\varphi(h)
\circ D(\lambda_{n+1}^+) (0)= D_\varphi(h) \circ D(\delta_n)\circ
D(\lambda_{n}^+) (0)$. Therefore, equations (\ref{la conju}) and
(\ref{cerca}) imply that $ \varphi \circ D(h\lambda_{n+1}^+) \circ
\varphi^{-1}(0)= D_\varphi(h\lambda_{n+1}^+) (0)$. Hence, equation
(\ref{distintos}) reads
$$D_\varphi(h\lambda_{n+1}^+)(0)>D_\varphi(g\lambda_{n+1}^+)(0),$$
which implies that $h\lambda_{n+1}^+ \succ^\prime g\lambda_{n+1}^+$.
In particular, we have that $\preceq^\prime$ is different from
$\preceq$ because we had assumed that $h\lambda_{n+1}^+\prec
g\lambda_{n+1}^+$. This finishes the proof of Theorem \ref{teo
finito generado}. $\hfill\square$

%%%%%%%%%%%%%%%%%%%%%%%%%%%%%%%%%%%%%%%%%%%%%%%%%%%%%%%%%%%%%%%%%%%%%

\subsection{The general case}
\label{sec extension por compacidad}

%%%%%%%%%%%%%%%%%%%%%%%%%%%%%%%%%%%%%%%%%%%%%%%%%%%%%%%%%%%%%%%%%%%%%%

\hspace{0.35 cm} There is a well-known criterion from Conrad-Ohnishi
\cite{conrad,ohnishi} stating that a group $\Gamma$ is
left-orderable if and only if for every finite family $f_1,\ldots,
f_k$, all of them different from the identity, there exist
$\eta_i\in \{-1,1\}$, $i=1,\ldots,k$, such that the identity is not
contained in the smallest semigroup containing
$\{f_1^{\eta_1},\ldots,f_k^{\eta_k}\}$. We will denote this
semigroup by $ \langle f_1^{\eta_1},\ldots,f_k^{\eta_k}\rangle^+$.

\vsp

In \cite[Proposition 1.4]{navas}, Navas shows that this criterion
(and the analogous one for bi-orderings \cite{ohnishi} and Conradian
orderings) is closely related to the compactness of
$\mathcal{LO}(\Gamma)$. Below, we present an extension of this
criterion that will permit us to deduce Theorem A from our proof of
Theorem \ref{teo finito generado}. This extension may be found in
\cite[Lemma 3.1.1]{kopytov}. However, for completeness, we give a
proof of it.

\vsp

Let $\gamma_1,\ldots, \gamma_n$ be a finite family of non-trivial
elements in a group $\Gamma$. We say that $\gamma_1,\ldots,
\gamma_n$ has property $(E)$ if and only if

\vsp

\noindent $(E)$: {\em for every finite family $f_1,\ldots, f_k$, of
elements different from the identity, there exists $\eta_i\in
\{-1,1\}$, $i=1,\ldots,k$, such that $id \not\in
\langle\gamma_1,\ldots,\gamma_n, f_1^{\eta_1},
\ldots,f_k^{\eta_k}\rangle^+$.}

\vsp

\noindent We say that such a choice of exponents $\eta_i$ is {\em
compatible}.

\begin{lem}\label{lema compacidad}{\em Let $\gamma_1,\ldots,\gamma_n$ be non trivial elements
in a left-orderable group $\Gamma$. Then $\Gamma$ admits a
left-ordering $\preceq$ such that $\gamma_i\succ id$ for all
$i=1,\ldots, n$, if and only if $\,\gamma_1,\ldots,\gamma_n$ has
property $(E)$. }
\end{lem}

\noindent {\em Proof:} The necessity of property $(E)$ is obvious.

\vsp

To see the sufficiency we will use the compactness of
$\mathcal{PLO}(\Gamma)$. For each finite family $f_1,\ldots,f_k$ of
non-trivial elements in $\Gamma$, and each compatible choice of
$\eta_i$, we let $\chi(f_1,\ldots,f_k;\eta_1,\ldots,\eta_k)$ be the
(closed) set of all partial-left-orderings such that each
$\gamma_j$, $j=1,\ldots,n$, and each $f_\ell^{\eta_\ell}$,
$\ell=1,\ldots,k$, is positive. By hypothesis, this set is
non-empty.

\vsp

Now, let $\chi(f_1,\ldots,f_k)$ be the (finite) union of all the
sets of the form $\chi(f_1,\ldots,f_k;\eta_1,\ldots,\eta_k)$, where
the choice of the exponents $\eta_i$ is compatible. Note that if $\{
\chi_i=\chi(f_{i,1},\ldots,f_{i,k}); \; i=1,\ldots n\}$ is a finite
family of subsets of this form, then, the intersection
$\chi_1\cap\ldots \cap \chi_n$ contains (the non-empty)
$\chi(f_{1,1},\ldots,f_{1,k}, \ldots, f_{n,1},\ldots,f_{n,k})$.
Since $\mathcal{PLO}(\Gamma)$ is compact, a direct application  of
the finite intersection property shows that $\chi$, the intersection
of all the sets of the form $\chi(f_{1},\ldots,f_{k})$, is
non-empty. It is quite clear that any partial-left-ordering
$\preceq\in \chi$ is a total ordering of $\Gamma$. Hence, any
left-ordering on $\chi$ is a left-ordering in which each $\gamma_i$,
$i=1,\ldots, n$, is positive $\hfill\square$

\vs

We now pass to the Proof of Theorem A.

\vs

Let $\preceq$ be a left-ordering on $G*H$, and let $F$ be a finite
subset of $\preceq$-positive elements in $G*H$ on which we want to
approximate $\preceq$. Let $G_0\subset G$ and $H_0\subset H$ be two
finite non-empty sets such that $G_0=G_0^{-1}$, $H_0=H_0^{-1}$ and
such that $F\subset \langle G_0\rangle * \langle H_0\rangle$. Let
$\Gamma_0=G_0\cup H_0=\Gamma_0^{-1}$ and $\Gamma=\langle
\Gamma_0\rangle$.

\vsp

Let $n\in \N$ be such that $F\subset B_{\Gamma_0}(n)$ and let $g\in
G_0$ and $h\in H_0$ be such that
$\lambda^+_{(B_{\Gamma_0}(n),\preceq)}\prec
h\lambda^+_{(B_{\Gamma_0}(n),\preceq)}$ and
$\lambda^+_{(B_{\Gamma_0}(n),\preceq)} \prec
\lambda^+_{(B_{\Gamma_0}(n),\preceq)}$ (see  Definition \ref{def
box}). As in the proof of Theorem \ref{teo finito generado}, we may
also assume that $h\lambda^+_{(B_{\Gamma_0}(n),\preceq)}\prec
g\lambda^+_{(B_{\Gamma_0}(n),\preceq)}$ (otherwise, we change the
names of $G$ and $H$). Finally, let
$\gamma_*=(h\lambda^+_{(B_{\Gamma_0}(n),\preceq)})^{-1}
g\lambda^+_{(B_{\Gamma_0}(n),\preceq)}$. Note that $id\prec
\gamma_*$.

% and that $\gamma_*\not\in B_{\Gamma_0}(n)$.

\vs

Theorem A follows directly from

\vs

\noindent {\bf Claim A:} The set $F\cup\{\gamma_*^{-1}\}$ has
property $(E)$.

\vs

In its turn, Claim A follows directly from

\begin{lem} \label{lema Claim A}{\em With the notations above, for any finitely generated subgroup
$\hat{\Gamma}$ of $G*H$ such that $\Gamma\subset \hat{\Gamma}$,
there exists a left-ordering $\preceq^*$ on $\hat{\Gamma}$ such that
any element in $F\cup\{\gamma_*^{-1}\}$ is $\preceq^*$-positive. }
\end{lem}

\noindent{\em Proof:} The proof follows the same lines as the proof
of Theorem \ref{teo finito generado}. Fix $\hat{\Gamma}$ a finitely
generated subgroup of $G*H$ containing $\Gamma$. We let
$\hat{\Gamma}_0$ be the generating set of $\hat{\Gamma}$. By
eventually enlarging $\hat{\Gamma}$, we shall assume that
$\hat{\Gamma}_0=\hat{G}_0\cup\hat{H}_0$, where $\hat{G}_0\subset G$
and $\hat{H}_0\subset H$, both non-empty sets. In this way we have
that $\hat{\Gamma}_0=\langle \hat{G}_0\rangle * \langle
\hat{H}_0\rangle$.

\vsp

To avoid heavy notation, for any $k\in \N$, we let $
\lambda_k^+=\lambda^+_{(B_{\Gamma_0}(k),\preceq)}$ and
$\lambda_k^-=\lambda^-_{(B_{\Gamma_0}(k),\preceq)}$.

\vsp

We let $D:\hat{\Gamma}\to Homeo_+(\R)$ be the dynamical realization
of the restriction of $\preceq$ to $\hat{\Gamma}$, that is, for any
$\gamma\in \hat{\Gamma}$, $\gamma\succ id$ if and only
$D(\gamma)(0)>0$.

\vsp

We let $x_0,x_1, y_0, y_1$ in $\R$ be such that
$$D(\lambda^+_{n+1})(0)<x_0<x_1<D(h\lambda_{n+1}^+)(0) <
D(g\lambda_{n+1}^+)(0)<y_1<y_0.$$

\vsp

We let $\varphi\in Homeo_+(\R)$ be such that $supp(\varphi)=\{x\in
\R\mid \varphi(x)\not=x\}=(x_0,y_0)$ and that $\varphi(x_1)>y_1$.
This implies that
\begin{equation}\label{distintos2}\varphi \circ D(h\lambda_{n+1}^+) \circ
\varphi^{-1}(0)> D(g\lambda_{n+1}^+)(0).\end{equation} Moreover, for
any $\bar{h}\in H_0$, and any $x\in
[D(\lambda_n^-)(0),D(\lambda_n^+)(0)]$, we have that
$D(\bar{h})(x)\leq D(\lambda_{n+1}^+)(0)<x_0$. Thus we conclude,
\begin{equation}\label{critico2}\varphi\circ D(\bar{h})\circ \varphi^{-1}(x)=D(\bar{h})(x),
 \text{ for all } x\leq D(\lambda_n^+)(0), \text{ and all } \bar{h}\in H_0.\end{equation}

\vsp

Now, let $D_\varphi:\hat{\Gamma} \to Homeo_+(\R)$ be defined by
$D_\varphi(\bar{g})=D(\bar{g})$ for all $\bar{g}\in
\langle\hat{G}_0\rangle$, and $D_\varphi(\bar{h})=\varphi\circ
D(\bar{h}) \circ \varphi^{-1}$ for all $\bar{h}\in \langle\hat{H}
\rangle$. Since $\hat{\Gamma}$ is the free product of
$\langle\hat{G}_0\rangle $ and $\langle\hat{H}_0\rangle$, we have
that $D_\varphi$ is an homomorphism (not necessarily injective).
Now, from the definition of $D_\varphi$ and equation
(\ref{critico2}), we have that $$ D(\gamma)(x)=D_\varphi(\gamma)(x),
\text{ for all $\gamma \in \Gamma_0$ and any $x\leq
D(\lambda_{n+1}^+)(0)$.}$$ In particular, for each $\gamma\in
\Gamma_0$, the graphs of $D(\gamma)$ and $D_\varphi(\gamma)$
coincide inside the square \break
$[D(\lambda_n^-)(0),D(\lambda_n^+)(0)]^2$. Hence, from Lemma
\ref{signo}, we have that
\begin{equation}\label{cerca2} \text{ for all } \gamma\in B_{\Gamma_0}(n),
\;\;D(\gamma)(0)=D_\varphi(\gamma)(0).\end{equation}

\vsp

Now, from Lemma \ref{lema extension}, there is a left-ordering
$\preceq^*$ on $\hat{\Gamma}$ such that $D_\varphi(\gamma)(0)
>0$ implies $\gamma\succ^* id$. Then, equation
(\ref{cerca2}) implies that $\preceq$ and $\preceq^*$ coincide on
$B_{\Gamma_0}(n)$. In particular, any element in $F$ is
$\preceq^*$-positive.

\vsp

However, arguing as in the end of the proof of Theorem \ref{teo
finito generado}, it can be shown that equation (\ref{distintos2})
is the same as
$$D_\varphi(h\lambda_{n+1}^+)(0)>D_\varphi(g\lambda_{n+1}^+)(0),$$
which shows that $id\prec^* \gamma_*^{-1}$.  $\hfill\square$

\vs

This finishes the proof of Theorem A.

%%%%%%%%%%%%%%%%%%%%%%%%%%%%%%%%%%%%%%%%%%%%%%%%%%%%%%%%%%%%%%%%%%%%%%%%%%%%%%%%%%%

\subsection{An example}
\label{sec ejemplo}

%%%%%%%%%%%%%%%%%%%%%%%%%%%%%%%%%%%%%%%%%%%%%%%%%%%%%%%%%%%%%%%%%%%%%%%%%%%%%%%%%%%%

\hspace{0.35 cm} We have proved that no left-ordering on a free
product of groups is isolated. In particular no positive cone of a
left-ordering on a free product is finitely generated as a semigroup
\cite[Proposition 1.8]{navas}. In this section, we show that there
exist a group with an isolated left-ordering whose positive cone is
not finitely generated as a semigroup. This seems to be the first
example of a group with this property.

\begin{prop} {\em The group $\Gamma=\langle a,b\mid bab^{-1}=a^{-2}\rangle$
is a finitely generated group with an isolated left-ordering whose
positive cone is not finitely generated as a semigroup.}
\end{prop}

\noindent {\em Proof:} The group $\Gamma$ is a group fitting in the
classification of groups having only finitely many left-orderings
\cite[Theorem 5.2.1]{kopytov}. However, we shall provide a direct
argument showing that it contains an isolated left-ordering.

\vsp

Let $\Gamma_1$ be the subgroup generated by $\{b^jab^{-j}\mid j\in
\Z\}$, and let $m, n$  in $\Z$. Note that both $b^n a b^{-n}$ and
$b^m ab^{-m}$ belong to $\langle b^k a b^{-k} \rangle$, where
$k=min\{0,n,m\}$. In particular, $\Gamma_1$ is an Abelian group
which is isomorphic to a non cyclic subgroup of the rational
numbers. Furthermore, $\Gamma_1$ is normal in $\Gamma$ and the
quotient $\Gamma/\Gamma_1=\langle b\Gamma_1 \rangle$ is isomorphic
to $\Z$.

\vsp

We let $\preceq_*$ be a left-ordering of $\Gamma_1$ such that
$a\succ_* id$, and $\preceq^*$ be a left-ordering on
$\Gamma/\Gamma_1$ such that $b\Gamma_1\succ^*\Gamma_1$. In this way
we can left-order $\Gamma$ by declaring
$$ g\succ id \Leftrightarrow \left\{
\begin{array}{l } g\Gamma_1\not=\Gamma_1 \text{ and } g\Gamma_1\succ^* \Gamma_1\;,  \text{ or }\\
g\in \Gamma_1 \;\;\text{ and } g\succ_* id .
\end{array} \right.$$
We claim that $\preceq$ is an isolated left-ordering. Indeed, let
$\preceq^\prime$ be a left-ordering such that $b\succ^\prime id$ and
such that $a\succ^\prime id$. In particular, since $\Gamma_1$ is
isomorphic to a subgroup of the rational numbers, we have that
$\preceq^\prime$ coincide with $\preceq$ on $\Gamma_1$. Now let
$g\in \Gamma$ be such that $g\notin \Gamma_1$. Let $n\in
\Z\setminus\{0\}$ be such that $b^n\Gamma_1=g\Gamma_1$, that is,
$g=b^n g_1$ for some $g_1\in \Gamma_1$. Suppose first that $n\geq
1$. In this case we have that $g=b^ng_1=b^{n-1}g_1^{-2}b$, which
shows that we can write $g$ as a product of $\preceq^\prime$-
positive elements. In particular $g\succ^\prime id$. In the case
that $n\leq -1$, the preceding argument shows that $g^{-1}$ is
$\preceq^\prime$-positive. Hence, we have that $\preceq^\prime$
coincide with $\preceq$ on $\Gamma$, showing that $\preceq$ is an
isolated left-ordering.

\vsp

Now, suppose by way of a contradiction that $\preceq$ has a positive
cone which is finitely generated as a semigroup. That is,
$P_\preceq=\{\gamma \in \Gamma \mid \gamma\succ id\}=\langle
S\rangle^+$, where $S=\{ \gamma_1,\ldots \gamma_n\}$. By eventually
re-labeling $S$, we may assume that $S=\{
\gamma_1,\ldots,\gamma_j,\ldots \gamma_n\}$, where $\gamma_i
\Gamma_1\succ^* \Gamma_1$, for $1\leq i\leq j$, and $\gamma_i \in
\Gamma_1$, for $i>j$. For $1\leq i \leq j$ we let
$\gamma_i=b^{n_i}g_i$, where $n_i\geq 1$ and $g_i\in \Gamma_1$.

\vsp

Now let $w=\gamma_{m_1}\ldots \gamma_{m_k}$ be an element in
$\langle S \rangle^+$. Since $\preceq^*$ is a left-ordering, we have
that $w\Gamma_1\succeq^* \Gamma_1$. This implies that any
$\preceq$-positive $g\in \Gamma_1$ may be written as a product of
$\gamma_{j+1},\ldots,\gamma_n$. However, this is impossible since
$\Gamma_1$ is a non-cyclic subgroup of the rational numbers. This
settles the desired contradiction. $\hfill\square$

%%%%%%%%%%%%%%%%%%%%%%%%%%%%%%%%%%%%%%%%%%%%%%%%%%%%%%%%%%%%%%%%%%%%%%%%%%%%%%%%%%%%%%%%%%%

\section{Constructing a dense orbit in the space of left-orderings of the free group}
\label{sec dense}

%%%%%%%%%%%%%%%%%%%%%%%%%%%%%%%%%%%%%%%%%%%%%%%%%%%%%%%%%%%%%%%%%%%%%%%%%%%%%%%%%%%%%%%%%%%%

\hspace{0.36 cm} We now proceed to the the construction of a
left-ordering on the free group of countable rank greater than one
whose orbit is dense under the natural conjugation action. The rough
idea is the following. Since the space of left-orderings of a
countable group is a compact metric space (see for instance
\cite{witte,navas,sikora} or the beginning of \S \ref{sec finito
generado}), it contains a dense countable subset. Now, we can
consider the {\em dynamical realization} (see \S \ref{real din}) of
each of these left-orderings, and cut large pieces from each one of
them (see for instance Definition \ref{def box}). Since we are
working with a free group, we can glue these pieces of dynamical
realizations together in a sole action of our group on the real
line. Moreover, if the gluing is made with a little bit of care,
then we can ensure very nice conjugacy properties from which we can
deduce Theorem B.

\vsp

First, we define an enumeration of the set of {\em balls} on a
countable free group. Let $S^+_\omega=\{a,b,\alpha_1,\alpha_2\ldots
\}$ be a free generating set of the free group of countable infinite
rank $F_\omega$. For $m\in\N=\{1,2,\ldots\}$, we let
$S_m^+=\{a,b,\alpha_1,\ldots \alpha_{m-2}\}$ if $m\geq 2$, and
$S^+_1=\{a\}$. For $n\in \N \cup \{\omega\}$, we let $S_n=S_n^+\cup
(S_n^+)^{-1}$. Note that we have the inclusion $S_n\subset S_w$, and
that $F_n=\langle S_n \rangle$. Using the notations of Definition
\ref{def box}, we let
$$ \mathcal{B}(F_n)=\left\{
\begin{array}{l } \{B_{S_n}(m)\mid m\in \N\} \;\text{ if } n\not=\omega\; , \\
\{ B_{S_m}(m)\mid m\in \N \}\; \text{ if } n=\omega.
\end{array} \right.$$  We call
$\mathcal{B}(F_n)$ the set of {\em balls} in $F_n$. We define
$\phi_n:\N\to \mathcal{B}(F_n)$ by $\phi_n(m)=B_{S_n}(m)$ if
$n\not=\omega$ and $\phi_\omega(m)=B_{S_{m+1}}(m+1)$. Note that,
$\cup_{m\in \N}\phi_n(m)=F_n$ and that, for any $B\in
\mathcal{B}(F_n)$, $n\not=1$, we have that $a$ and $b$ belong to
$B$. Note also that $S_\omega\cap \phi_n(m)=S_k$, where
$k=min\{n,m+1\}$ (assuming that $\omega$ is bigger than any
integer).

\vs

Fix once and for all $n\in \N \cup\{\omega\}$, $n\not=1$. Let
$\phi=\phi_n$, $\mathcal{B}=\mathcal{B}(F_n)$, and
$\mathcal{D}=\{\preceq_1,\preceq_2,\ldots\}$ be a countable dense
subset of $\mathcal{LO}(F_n)$. Let $\eta:\Z\to \mathcal{B}\times
\mathcal{D}$ be a surjection, with
$\eta(k)=(\phi(n_k),\preceq_{m_k})$.

\vsp

By Remak \ref{rem cajas} we have that there exists
$D_{\eta(k)}:F_n\to Homeo_+(\R)$, a dynamical realization-like
homomorphism for $\preceq_{m_k}$, such that:

\vsp

\noindent $(i)$ The reference point for $D_{\eta(k)}$ is $k$.

\vsp

\noindent $(ii)$ The $\eta(k)$-box coincides with the square
$[k-1/3, k+1/3]^2$.

\vs

Theorem B is a direct consequence of the following

\begin{prop} \label{the action} {\em There is an homomorphism $D:F_n\to Homeo_+(\R)$ such
that, for each $k\in \Z$, inside $[k-1/3,k+1/3]^2$, the graphs of
$D(g)$ coincide with the graphs of $D_{\eta(k)}(g)$ for any $g\in
S_n\cap\phi(n_k)$. In this action, all the integers lie in the same
orbit.}
\end{prop}

\noindent {\em Proof of Theorem B from Proposition \ref{the
action}:} Let $(x_0,x_1,\ldots )$ be a dense sequence in $\R$ such
that $x_0=0$ (note that $0$ may not have a free orbit), and let $D$
be the homomorphism given by Proposition \ref{the action}. Note that
$D$ is an embedding, since, from Lemma \ref{signo}, we have that any
non-trivial $w\in \phi(n_k)$ acts nontrivially at the point $k\in
\R$. Hence, we may let $\preceq$ be the induced left-ordering on
$F_n$ from the action $D$ and the reference points $(
x_0,x_1,x_2,\ldots)$. In particular, for $h\in F_n$, we have that
$D(h)(0)>0\Rightarrow h\succ id$. We claim that $\preceq$ has a
dense orbit under the natural action of $F_n$ on
$\mathcal{LO}(F_n)$.

\vsp

Clearly, to prove our claim it is enough to prove that the orbit of
$\preceq$ accumulates at every $\preceq_m \in \mathcal{D}$. That is,
given $\preceq_m$ and any finite set $\{ h_1,h_2,... ,h_N \}$ such
that $id \prec_m h_j$, for $1\leq j \leq N$, we need to find  $w\in
F_n$ such that $h_j\succ_w id$ for every $1\leq j\leq N$, where, as
defined in the Introduction, $h\succ_w id $ if and only if
$whw^{-1}\succ id$.

\vsp

Let $j\in \N$ be such that $h_1,\ldots , h_N$ belongs to $\phi(j)$.
Let $k$ be such that $\eta(k)=(\phi(j),\preceq_m)$. By Proposition
\ref{the action}, there is $w_k\in F_n$ such that $D(w_k)(0)=k$.
Also by Proposition \ref{the action}, inside $[k-1/3,k+1/3]^2$, for
every $g\in S_\omega\cap\phi(j)$ we have that the graphs of $D(g)$
are the same as those of $D_{\eta(k)}(g)$. Then, Lemma \ref{signo}
implies that for each $h_j$, $1\leq j\leq N$, we have that
$h_i\succ_m id $ if and only if $D(h_j)(k)>k$. But this is the same
as saying that $D(h_j) (D(w_k)(0))> D(w_k)(0)$, which implies that
$D(w_k^{-1})\circ D(h_j)\circ D(w_k) (0)>0$. Therefore, by
definition of $\preceq$, we have that $w_k^{-1} h_j w_k \succ id$
for every $1\leq j \leq N$. Now, by definition of the action of
$F_n$ on $\mathcal{LO}(F_n)$, this implies that $\preceq_{w^{-1}_k}$
is a left-ordering such that $h_j\succ_{w^{-1}_k} id$. This finishes
the proof of Theorem B. $\hfill\square$

\vsp\vsp

To prove Proposition \ref{the action} we first consider $g\in S_n$,
and let $K=\{k\in \Z \mid g\in \phi(n_k)\}$. Now if $k_0$ and $k_1$
are elements of $K$ such that $k_0<k_1$ and such that there is no
other element of $K$ in between, then we can linearly interpolate
the portion of the graph of $D_{\eta(k_0)}(g)$ inside
$[k_0-1/3,k_0+1/3]^2$ until the portion of the graph of
$D_{\eta(k_1)}(g)$ inside $[k_1-1/3,k_1+1/3]^2$. Repeating this
argument, we get a function $\hat{g}\in Homeo_+(\R)$ that coincides
with $D_{\eta(k)}(g)$ for all $k\in K$. In this way we have proved

\begin{lem} \label{inductivo 2}{\em Let $g\in S_n$. For each $k\in \Z$ we let $n_k$ and $m_k$ in $\N$ be
such that $\eta(k)=(\phi(n_k),\preceq_{m_k})$. Then, there exist
$\hat{g}\in Homeo_+(\R)$ such that for every $k\in \Z$ such that
$g\in \phi(n_k)$, the graph of $\hat{g}$ inside $[k-1/3,k+1/3]^2$
coincide with the graphs of $D_{\eta(k)}(g)$.}
\end{lem}

\begin{lem} \label{inductivo} {\em For each $k\in \Z$, we can modify the
homeomorphisms $\hat{a}$ and $\hat{b}$ (given by Lemma
\ref{inductivo 2}) inside $[k-1/3, k+1+1/3]^2$ but outside
$[k-1/3,k+1/3]^2 \cup [k+1-1/3,k+1+1/3]^2$ (see Figure 4.1) in such
a way that the modified homeomorphisms, which we still denote
$\hat{a}$ and $\hat{b}$, have the following property

\vsp\vsp

\noindent $(P):$ there is a reduced word $w$ in the free group
generated by $\{ \hat{a}, \hat{b}\}$ such that $w(k+1/3)=k+1-1/3$.
Moreover, the iterates of $k+1/3$ along the initial segments of $w$
remain inside $[k-1/3,k+1+1/3]$.}
\end{lem}

\vspace{0.8cm}

%%%%%%%%%%%%%%%%%%%%%%%%%%%%%%%%%%%%%%%%%%%%%%%%%%%%%%%%%%%%%%%%%%%%%%%%%%%%%%%%%%%%%%%%%%%%

\beginpicture

\setcoordinatesystem units <1cm,1cm>

%%%%%%%%%%%%%%%%%%%%%%%%%%%%%%%%%%%%%%%%%%%%%%%%%%%%%%%%%%%%%%%%%%%%%%%%%%%%%%%%%%%%%%%%%%%%%%%

\putrule from 2.5 0 to 4 0 \putrule from 4 0 to 4 -1.5 \putrule from
2.5 -1.5 to 4 -1.5 \putrule from 2.5 -1.5 to 2.5 0

\putrule from 4.5 0.5 to 6 0.5 \putrule from 4.5 0.5 to 4.5 2
\putrule from 6 0.5 to 6 2 \putrule from 4.5 2 to 6 2

\putrule from 0.5 -2.1 to 8 -2.1

%%%%%%%%%%%%%%%%%%%%%%%%%%%%%%%%%%%%%%%%%%%%%%%%%%%%%%%%%%%%%%%%%%%%%%%%%%%%%%%%%%%%%%%%%%%%%%
\setdots

\plot 2.5 -1.5 2.5 2 / \plot 2.5 2  6 2 / \plot 6 2  6 -1.5 / \plot
6 -1.5 2.5 -1.5 /

\plot 1.8 -2.2 6.4 2.4 /

\plot 4 0  4 -2.8 / \plot 4 0  1.4  0 /

\plot 6 0.5 7 0.5 / \plot 4.5 0.5 4.5 2.5 /

\put{Figure 4.1} at 4 -3.5 \put{} at -4.2 0

\small

\put{$\eta(k)$} at 3.1 -0.7 \put{-box} at 3.4 -1.1

\put{$\eta(k+1)$} at 5.3 1.3 \put{-box} at 5.3 1

\put{$\bullet$} at  5.3 -2.1 \put{$\bullet$} at  3.3 -2.1
\put{$k+1$} at  5.3 -2.4 \put{$k$} at  3.3 -2.4

\put{m} at 2.9 1.2 \put{o} at 3.2 1.05 \put{d} at 3.4 0.9 \put{i} at
3.7 0.7 \put{f} at 4 0.5 \put{i} at 4.25 0.3 \put{c} at 4.5 0.1
\put{a} at 4.75 -0.1 \put{t} at 5 -0.2 \put{i} at 5.2 -0.3 \put{o}
at 5.4 -0.4 \put{n} at 5.6 -0.55 \put{s} at 5.8 -0.7

\endpicture

%%%%%%%%%%%%%%%%%%%%%%%%%%%%%%%%%%%%%%%%%%%%%%%%%%%%%%%%%%%%%%%%%%%%%%%%%%%%%%%%%%%%%%%%%%%%%%%%
%%%%%%%%%%%%%%%%%%%%%%%%%%%%%%%%%%%%%%%%%%%%%%%%%%%%%%%%%%%%%%%%%%%%%%%%%%%%%%%%%%%%%%%%%%%%%%%%

\vspace{0.8cm}

\noindent {\em Proof:} For $h\in \{\hat{a}^{\pm 1},\, \hat{b}^{\pm
1}\}$, define $l_h=\sup\{ x\in [k-1/3,k+1/3] \mid h(x)\leq k+1/3 \}$
and $r_h= \inf \{ x \in [k+1-1/3,k+1+1/3] \mid h(x)\geq k+1-1/3\}$.
Let $x_0\in\, ]k+1/3, k+1-1/3[$. To modify $\hat{a}$ and $\hat{b}$,
we proceed as follows:

\vsp

\noindent \textbf{Case 1:} There is $h\in \{\hat{a}^{\pm1},\,
\hat{b}^{\pm1}\}$ such that $l_h< k+1/3$ and $r_h=k+1-1/3$.

\vsp

In this case, we (re)define $h$ linearly from $(l_h, h(l_h))=(l_h,
k+1/3)$ to $( k+1/3 ,x_0)$, then linearly from $( k+1/3 ,x_0)$ to
$(x_0, k+1-1/3)$, and then linearly from $(x_0,k+1-1/3)$ to
$(k+1-1/3, h(k+1-1/3))= (r_h, h(r_h))$; see Figure 4.2 (a). The
other generator, say $f$, may be extended linearly from
$(l_f,f(l_f))$ to $(r_f, f(r_f))$.

\vsp

Note that in this case we have $h(k+1/3)=x_0$ and $h(x_0)=k+1-1/3$.
This shows that $(P)$ holds for $w=h^2$.

\vs

We note that, for $h\in \{\hat{a}^{\pm1},\, \hat{b}^{\pm1}\}$, we
have that $l_h=k+1/3 \Leftrightarrow l_{h^{-1}}<k+1/3\;$ and $\;
r_h=k+1-1/3\Leftrightarrow r_{h^{-1}}>k+1-1/3$. Therefore, if there
is no $h$ as in Case 1, then we are in

%Also note that we have the liberty to extend the second generator in
%any way we want. This would become important for prove that the
%ordering whose orbit is dense has a convex subgroup.

\vsp

\noindent \textbf{Case 2:} There are  $f,h \in \{\hat{a}^{\pm1},
\hat{b}^{\pm 1} \}$ such that $l_h< k+1/3$, $\;r_h>k+1-1/3$,
$\;l_f<k+1/3$ and $r_f>k+1-1/3$.

\vsp

In this case we define $h$ linearly from $(l_h,h(l_h))$ to $(k+1/3,
x_0)$, and then linearly from $(k+1/3, x_0)$ to $(r_h, h(r_h))$. For
$f$, we define it linearly from $(l_f, f(l_f))$ to $(k+1-1/3, x_0)$,
and then linearly from $(k+1-1/3, x_0)$ to $(r_f,f(r_f))$; see
Figure 4.2 (b).

\vsp

Note that $h(k+1/3)=x_0=f(k+1-1/3)$. This shows that $(P)$ holds for
$w=f^{-1}h$. $\hfill\square$

\vspace{0.8cm}

%%%%%%%%%%%%%%%%%%%%%%%%%%%%%%%%%%%%%%%%%%%%%%%%%%%%%%%%%%%%%%%%%%%%%%%%%%%%%%%%%%%%%%%%%%%%

\beginpicture

\setcoordinatesystem units <1cm,1cm>

%%%%%%%%%%%%%%%%%%%%%%%%%%%%%%%%%%%%%%%%%%%%%%%%%%%%%%%%%%%%%%%%%%%%%%%%%%%%%%%%%%%%%%%%%%%%%%%
\plot 1.42 1.7 1.8 2.1 / \plot 1.8 2.1 2.4 2.2 / \plot -1.7 -2.4
-1.3 -1.4 / \plot -1.3 -1.4 -1 0 / \plot -1 0  0.2 1.35 / \plot 0.2
1.35 1.42 1.7 /

\plot 7.5 -2 8 -1.24 / \plot 8 -1.24 8.55 0 / \plot 8.55 0 11.5 1.27
/ \plot 11.5 1.27 12.3 2 /

\plot 7.9 -2.2 8.3 -1.24 / \plot 8.3 -1.24 10.95 0 / \plot 10.95 0
11.2 2.3 /
%%%%%%%%%%%%%%%%%%%%%%%%%%%%%%%%%%%%%%%%%%%%%%%%%%%%%%%%%%%%%%%%%%%%%%%%%%%%%%%%%%%%%%%%%%%
\putrule from -2.5 0 to 2.8 0 \putrule from 0.2 2.5 to 0.2 0
\putrule from 0.2 -2.7 to 0.2 -0.5

\putrule from -2.5 -1.35 to -1 -1.35 \putrule from -1 -1.35 to -1
-2.7

\putrule from 1.4 2.7 to 1.4 1.35 \putrule from 1.4 1.35 to 2.9 1.35

\putrule from 7 0 to 12.5 0 \putrule from 9.74 2.5 to 9.74 0
\putrule from 9.74 -2.7 to 9.74 -0.5

\putrule from 7 -1.24 to 8.55 -1.24 \putrule from 8.55 -1.24 to 8.55
-2.7

\putrule from 10.95 2.7 to 10.95 1.27 \putrule from 10.95 1.27 to
12.5 1.27

%%%%%%%%%%%%%%%%%%%%%%%%%%%%%%%%%%%%%%%%%%%%%%%%%%%%%%%%%%%%%%%%%%%%%%%%%%%%%%%%%%%%%%%%%%%%%%
\setdots \plot -2.3 -2.7 2.7 2.7 /    \plot 7 -2.8 12.5 2.8 /

\put{Figure 4.2 (a)} at 0.3 -3.5 \put{Figure 4.2 (b)} at 9.8 -3.5
%%%%%%%%%%%%%%%%%%%%%%%%%%%%%%%%%%%%%%%%%%%%%%%%%%%%%%%%%%%%%%%%%%%%%%%%%%%%%%%%%%%%%%%%%%%%%%%%
\small

\put{$x_0$} at 0.2 -0.3 \put{$\bullet$} at 0.2  0

\put{$l_h$} at -1.45 -0.3 \put{$\bullet$} at -1.3 0
%\put{$(l_h,h(l_h))$} at -1.95 -1 \put{$\bullet$} at -1.3 -1.37
\put{$k+\frac{1}{3}$} at -1.2 0.5 \put{$\bullet$} at -1 0
\put{$k+1-\frac{1}{3}$} at 1.6 -0.3 \put{$\bullet$} at 1.4 0

\put{$x_0$} at 9.78 -0.3 \put{$\bullet$} at 9.75  0 \put{$\bullet$}
at 8.05 -0.01 \put{$l_h$} at 7.9 -0.3 \put{$\bullet$} at 11.45 0
\put{$r_h$} at 11.45 0.3

\put{$k+\frac{1}{3}$} at 8.6 0.5 \put{$\bullet$} at 8.55 0
\put{$k+1-\frac{1}{3}$} at 11.3 -0.5 \put{$\bullet$} at 10.95 0

\put{$h$} at -1.5 -2.6 \put{$h$} at 2.6 2.2

\put{$h$} at 7.3 -2.1  \put{$f$} at 7.9 -2.5 \put{$h$} at 12.5 2.1
\put{$f$} at 11.25 2.6

\put{} at -3 0
\endpicture

%%%%%%%%%%%%%%%%%%%%%%%%%%%%%%%%%%%%%%%%%%%%%%%%%%%%%%%%%%%%%%%%%%%%%%%%%%%%%%%%%%%%%%%%%%%%%%%%
%%%%%%%%%%%%%%%%%%%%%%%%%%%%%%%%%%%%%%%%%%%%%%%%%%%%%%%%%%%%%%%%%%%%%%%%%%%%%%%%%%%%%%%%%%%%%%%%

\vspace{0.5cm}

\noindent {\em Proof of Proposition \ref{the action}:} For each
$g\in S^+_n$, we let $\hat{g}$ be as in Lemma \ref{inductivo 2}.
Hence, inside $[k-1/3,k+1/3]^2$, the graphs of $\hat{g}$ coincide
with the graphs of $D_{\eta(k)}(g)$ for any $g\in S_n\cap\phi(n_k)$,
where $\eta(k)=(\phi(n_k),\preceq_{m_k})$. Now, for each $k\in \Z$
we apply inductively Lemma \ref{inductivo} to modify $\hat{a}$ and
$\hat{b}$. This modified homeomorphisms will be still denoted
$\hat{a}$ and $\hat{b}$. Note that Lemma \ref{inductivo} implies
that the modifications are made in such a way that they do not
overlap one with each other and that, for each $k\in \Z$, the graphs
of $\hat{a}$ and $\hat{b}$ coincides with the graphs of
$D_{\eta(k)}(a)$ and $D_{\eta(k)}(b)$ inside $[k-1/3,k+1/3]^2$.
Therefore, if we define $\hat{D}:F_n\to Homeo_+(\R)$ by
$\hat{D}(g)=\hat{g}$ for every $g\in S_n$, we have that, inside
$[k-1/3,k+1/3]^2$, the graphs of $\hat{D}(g)$ coincide with the
graphs of $D_{\eta(k)}(g)$ for any $g\in S_n\cap\phi(n_k)$.

\vsp

Finally, for each $k\in \Z$, Lemma \ref{signo} implies that in the
action given by $\hat{D}$, the points $k,k+1/3$ and $k-1/3$ are in
the same orbit. Hence, from Lemma \ref{inductivo}, we have that in
this action all the integers are in the same orbit. This finishes
the proof of Proposition \ref{the action}. $\hfill\square$

%%%%%%%%%%%%%%%%%%%%%%%%%%%%%%%%%%%%%%%%%%%%%%%%%%%%%%%%%%%%%%%%%%%%%%%%%%%%%%%%%%%%%%%%%%%%%%%%%%%%%%%%%
%%%%%%%%%%%%%%%%%%%%%%%%%%%%%%%%%%%%%%%%%%%%%%%%%%%%%%%%%%%%%%%%%%%%%%%%%%%%%%%%%%%%%%%%%%%%%%%%%%%%%%%%%

\begin{small}

\vspace{0.1cm}

%%%%%%%%%%%%%%%%%%%%%%%%%%%%%%%%%%%%%%%%%%%%%%%%%%%%%%%%%%%%%%%%%%%%%%%%%%%%%%%%%%%%%%%%%

\vspace{0.37cm}

\noindent Crist\'obal Rivas\\

\noindent Dep. de Matem\'atica y C.C., Fac. de Ciencia, Univ. de Santiago de Chile\\

\noindent Alameda 3363, Estaci\'on Central, Santiago, Chile\\

\noindent Email: cristobalrivas@u.uchile.cl

\end{small}

%%%%%%%%%%%%%%%%%%%%%%%%%%%%%%%%%%%%%%%%%%%%%%%%%%%%%%%%%%%%%%%%%%%%%%%%%%%%%%%%%%%%%%%%%%
%%%%%%%%%%%%%%%%%%%%%%%%%%%%%%%%%%%%%%%%%%%%%%%%%%%%%%%%%%%%%%%%%%%%%%%%%%%%%%%%%%%%%%%%%%


\begin{thebibliography}{Dillo 83}

\bibitem{clay 2} {\sc A. Clay.} Free lattice ordered groups and the topology on
the space of left orderings of a group. To appear in {\em
Monatshefte für Mathematik} (available online).

\bibitem{conrad2} {\sc P. Conrad.} Free lattice-ordered groups. {\em J.
Algebra} {\bf 16} (1970), 191-203.

\bibitem{conrad} {\sc P. Conrad.} Right-ordered groups. {\em Mich. Math. Journal}
{\bf 6} (1959), 267-275.

\bibitem{braids} {\sc P. Dehornoy, I. Dynnikov, D. Rolfsen \& B.
Wiest.} \textit{Ordering Braids}. Math. Surveys and Monographs
\textbf{148}, A.M.S. (2008).

\bibitem{dd} {\sc T.V. Dubrovina \& N.I. Dubrovin}. On braid
groups. {\em Mat. Sb.} {\bf 192} (2001), 693-703.

\bibitem{ghys} {\sc \'E. Ghys.} Groups acting on the circle.
\textit{L'Enseignement Math\'ematique} {\bf 47} (2001), 239-407.

\bibitem{hockin young} {\sc J.G. Hocking \& G.S. Young.} {\em Topology.}
Addison-Wesley Publishing Co., Inc., Reading, Mass.-London (1961).

\bibitem{ito} {\sc T. Ito.} Dehornoy-like left orderings and isolated left orderings,
{\em Preprint.} arXiv 1102.4669v1.

\bibitem{kopytov2} {\sc V. Kopytov.} Free lattice-ordered groups.
\textit{Sibirsk Mat. Zh.} {\bf 24}(1) (1983), 120-124.

\bibitem{kopytov} {\sc V. Kopytov \& N. Medvedev.} {\em Right ordered groups.}
Siberian School of Algebra and Logic, Plenum Publ. Corp., New York
(1996).

\bibitem{linnell} {\sc P. Linnell.} The space of left orders of a group is either
finite or uncountable. {\em London Math. Soc.} {\bf 43} (2011),
200-202.

\bibitem{mccleary} {\sc S.H. McCleary.} Free lattice-ordered group represented as $o$-2 transitive $l$-permutation groups. \textit{Trans. Amer. Math. Soc.}, {\bf 290(1)} (1985),
69-79.

\bibitem{witte} {\sc D. Morris-Witte.} Amenable groups that act on the line.
{\em Algebr. Geom. Topol.} {\bf 6} (2006), 2509-2518.

\bibitem{navas-hecke} {\sc A. Navas.} A remarkable family of left-ordered
groups: central extensions of Hecke groups. {\em J. Algebra} {\bf
328} (2011), 31-42.

\bibitem{navas} {\sc A. Navas.} On the dynamics of (left) orderable
groups. {\em Ann. Inst. Fourier (Grenoble)} {\bf 60} (2010),
1685-1740.

\bibitem{navas bertol} {\sc A. Navas \& B. Wiest.} Nielsen-Thurston orderings and the space of
briad orderings. To appear in {\em Bull. of Lon. Math. Soc.}
(available online).

\bibitem{ohnishi} {\sc M. Ohnishi.} Linear order on a group. {\em
Osaka Math. J.} {\bf 2} (1959), 17-18.

\bibitem{rivas2} {\sc C. Rivas.} On groups with finitely many Conradian orderings. To appear in {\em
Comm. in Algebra}.

\bibitem{tesis} {\sc C. Rivas.} Orderable groups. Ph.D. thesis,
Univ. de Chile (2010).

\bibitem{sikora} {\sc A. Sikora.} Topology on the spaces of orderings of groups.
{\em Bull. London Math. Soc.} {\bf 36} (2004), 519-526.


\end{thebibliography}
\end{document}